\documentclass[letterpaper, 10 pt, conference,twocolumn]{ieeeconf}

\IEEEoverridecommandlockouts 

\usepackage{amsmath,amsfonts,amssymb}
\usepackage{amsthm}
\usepackage{graphicx,epsfig,psfrag,remark}
\usepackage{xfrac,bm}
\usepackage{algorithmic,algorithm}
\usepackage[all]{xy}
\usepackage{varioref}
\usepackage{wrapfig}
\usepackage{threeparttable}
\usepackage{dcolumn}
\newcolumntype{d}{D{.}{.}{-1}}
\usepackage{nomencl}
\makeglossary
\usepackage{caption}
\usepackage{subcaption}
\usepackage{comment,cite}
\usepackage{xfrac}
\usepackage{mathrsfs}
\usepackage{xcolor}

\newcommand{\x}{{\bm{x}}}
\newcommand{\X}{{\mathsf{X}}}

\newcommand{\f}{{\mathsf{f}}}

\newcommand{\z}{{\bm{z}}}

\newcommand{\cC}{\mathcal{C}}

\newcommand{\cJ}{\mathcal{J}}
\newcommand{\cK}{\mathcal{K}}
\newcommand{\cL}{\mathcal{L}}

\newcommand{\cN}{\mathcal{N}}

\newcommand{\cV}{\mathcal{V}}

\newcommand{\tr}[1]{\ensuremath{\operatorname{tr}\left( #1 \right)}}
\renewcommand{\d}{\ensuremath{\operatorname{d}}}

\theoremstyle{plain}
\newtheorem{proposition}{Proposition}

\newtheorem{problem}{Problem}
\newtheorem{theorem}{Theorem}

\renewcommand{\t}{\ensuremath{^{\mathrm{T}}}}

\newremark{remark}{Remark}

\graphicspath{{./plots}}

\renewcommand{\baselinestretch}{0.96}

\title{Covariance Control of Discrete-Time Gaussian Linear Systems\\ Using Affine Disturbance Feedback Control Policies}

\author{Isin M. Balci \thanks{I. M. Balci is a PhD student in the Department of Aerospace Engineering and Engineering Mechanics, The University of Texas at Austin, Austin. Texas 78712-1221, USA, Email: isinmertbalci@utexas.edu} \and Efstathios Bakolas \thanks{E. Bakolas
is an Associate Professor in the Department of Aerospace Engineering
and Engineering Mechanics, The University of Texas at Austin,
Austin, Texas 78712-1221, USA, Email: bakolas@austin.utexas.edu}}

\begin{document}
\maketitle

\begin{abstract}
    In this paper, we present a new control policy parametrization for the finite-horizon covariance steering problem for discrete-time Gaussian linear systems (DTGLS) which can reduce the latter stochastic optimal control problem to a tractable optimization problem. The covariance steering problem seeks for a feedback control policy that will steer the state covariance of a DTGLS to a desired positive definite matrix in finite time.
    We consider two different formulations of the covariance steering problem, one with hard terminal LMI constraints (hard-constrained covariance steering) and another one with soft terminal constraints in the form of a terminal cost which corresponds to the squared Wasserstein distance between the actual terminal state (Gaussian) distribution and the desired one (soft-constrained covariance steering).
    We propose a solution approach that relies on the affine disturbance feedback parametrization for both problem formulations. 
    We show that this particular parametrization allows us to reduce the hard-constrained covariance steering problem into a  semi-definite program (SDP) and the soft-constrained covariance steering problem into a difference of convex functions program (DCP).
    Finally, we show the advantages of our approach over other covariance steering algorithms in terms of computational complexity and computation time by means of theoretical analysis and numerical simulations.
\end{abstract}

\section{Introduction}\label{s:intro}
In this work, we consider the problem of characterizing computationally tractable control policies that will steer the mean and covariance of the terminal state of a discrete-time linear stochastic system ``close'' to respective goal quantities. This type of problems are referred to as \textit{covariance steering} (or \textit{covariance control}) problems in the literature of stochastic control.
We will consider two variations of the covariance steering problem. 
The goal in the first problem formulation is to steer the mean of the terminal state to a prescribed vector and have the terminal state covariance satisfy a certain LMI-type constraint; we refer to this problem as the hard constrained covariance steering (HCCS) problem.
In the second formulation of the covariance steering problem, we seek for a control policy that will minimize the distance between the terminal (Gaussian) distribution of the state and a desired goal (Gaussian) distribution measured in terms of the (squared) Wasserstein distance between the two distributions while satisfying the probabilistic input and state constraints; we refer to this problem as the soft-constrained covariance steering (SCCS) problem. 

\textit{Literature Review:} Infinite-horizon covariance control problems for both
continuous-time and discrete-time stochastic linear systems have been addressed in 
\cite{p:skeltonIJC,p:skeltonTAC,p:levy1997discrete}.
Finite-horizon covariance steering problems have recently received significant attention for both the continuous-time case 
\cite{p:georgiou15A,p:georgiou15B} and the discrete-time case
~\cite{p:bakcdc16,p:PT2017,p:BAKOLAS2018}.
Covariance control problems for the partial information case have been studied in 
\cite{p:EBACC17,p:bakTAC2019,p:ridderhof2020,p:kotsalis2020}.
Approaches that consider soft-constrained formulations of the covariance steering problem (based on appropriate terminal costs) can be found in~\cite{p:halderacc16,p:balci2020} and~\cite{p:grune2019}, in which the terminal cost is defined in terms of, respectively, the squared Wasserstein distance and the squared $\mathcal{L}_2$ spatial norm between the goal distribution and the distribution attained by the terminal state.

In our previous work, we have addressed covariance steering problems for discrete-time stochastic (Gaussian) linear systems under both full state and partial state information based on techniques from convex optimization~\cite{p:bakcdc16,p:EBACC17,p:BAKOLAS2018,p:bakTAC2019} and difference of convex functions programming~\cite{p:balci2020}.
In these references, the reduction of the stochastic optimal control problems to (finite-dimensional) optimization problems relied on the utilization of the state history feedback control parametrization~\cite{p:GOULART2006}. 
In this control policy parametrization, the control input at each stage is expressed as an affine function of the complete history of states visited by the system (including the current state). By using the state feedback control policy parametrization, one can reduce the covariance steering problem into a convex program via a bilinear transformation.~\cite{p:skaf2010}. Because, the whole history of states is used in this parametrization, the dimension of the resulting optimization problem can be prohibitive for problems with long (yet finite) time horizons.


\textit{Main Contribution:} In this paper, we present a new solution approach to the covariance steering 
problem (for both the hard-constrained and the soft-constrained problem formulations) in the case of full state information. Our approach is based on a control policy parameterization which can be interpreted as a modified version of the stochastic version of the affine disturbance feedback control parametrization~\cite{p:bennemiro2004} tailored to the covariance steering problem. We show that by using this particular control policy
parametrization, one can directly reduce the HCCS problem into a convex optimization problem and the SCCS problem into a difference of convex functions program, whose decision variables are essentially the controller parameters. 
This is in sharp contrast with the state history feedback control parametrization which requires significant pre-processing in order to associate the controller parameters (decision variables of the stochastic optimal control problem) with the decisions variables of the corresponding (finite-dimensional) optimization program by means of bilinear transformations. 
%
%
The fact that in our proposed approach the decision variables of the control and optimization problems are in direct correspondence also allows us to consider modified control policies which are based on truncated histories of the disturbances which have acted upon the system. Using these modified policies lead to more computationally tractable optimization problems (problems with fewer decision variables).

\textit{Structure of the paper:} 
The rest of the paper is organized as follows. 
In Section \ref{s:problem-formulation}, we formulate the two variations of the covariance steering problem and introduce our proposed policy parametrization. 
The reduction of the stochastic optimal control problem into a convex semidefinite program (for the HCCS problem) and a difference of convex functions program (for the SCCS problem) are described in Sections \ref{s:reduction2convex} and \ref{s:reduction2dc}, respectively. 
In Section \ref{s:experiments}, we present numerical simulations and finally, we conclude the paper in Section \ref{s:conclusion}.


\section{Problem Formulation}\label{s:problem-formulation}
\subsection{Notation}
We denote by $\mathbb{R}^n$ the set of $n$-dimensional real vectors
and by $\mathbb{Z}$ the set of integers. 
We write $\mathbb{E}\left[
\cdot\right]$ to denote the expectation functional. We denote the probability of the random event $E$ as $\mathbb{P}(E)$. Given two integers $\tau_1, \tau_2$ with $\tau_1\leq \tau_2$, then $[\tau_1, \tau_2]_d := [\tau_1,\tau_2] \cap \mathbb{Z}$. 
Given a finite sequence of vectors $\mathscr{X}:=\{x_i:~i\in[1,m]_d\}$, we denote by $\mathrm{vertcat}(\mathscr{X})$ the concatenation of its vectors, that is, $\mathrm{vertcat}(\mathscr{X}) := [x_1\t, \dots, x_m\t]\t$. 
Given a square matrix $\mathbf{A}$, we denote its trace by
$\operatorname{tr}(\mathbf{A})$, its Frobenius norm by $\lVert\mathbf{A} \rVert_{F}$, and its nuclear norm by $\lVert \mathbf{A} \rVert_{*}$. 
We write $\mathbf{0}$ and $I_n$ to denote the zero matrix (of suitable dimensions) and the $n\times n$ identity matrix, respectively. 
The space of real symmetric $n\times n$ matrices will be denoted by $\mathbb{S}_n$. Furthermore, we will denote the convex cone of
$n\times n$ (symmetric) positive semi-definite and (symmetric) positive definite matrices by $\mathbb{S}^{+}_n$ and
$\mathbb{S}^{++}_n$, respectively. 
Given $A \in \mathbb{S}^+_n$, we denote by $A^{1/2} \in \mathbb{S}^{+}_n$ its (unique) square root, that is, $A^{1/2} A^{1/2} = A$. 
We write $\mathrm{bdiag}(A_1,$ $\dots, A_\ell)$ to denote the block diagonal matrix formed by the matrices $A_i$, $i\in \{1,\dots, \ell\}$. 
Finally, we denote by $\mu_z$ and $\mathrm{var}_z$ the mean and the covariance (or variance) of a random vector $z$, respectively, that is, $\mu_z:=\mathbb{E}[z]$ and $\mathrm{var}_z:= \mathbb{E}[(z-\mu_z)(z-\mu_z)\t]=\mathbb{E}[zz\t]-\mu_z \mu_z\t$.

\subsection{Squared Wasserstein Distance}
The Wasserstein distance defines a valid metric (i.e., satisfies all relevant axioms of a metric) in the space of probability distributions. Although, in general, it is not possible to find a closed-form expression for the Wasserstein distance between two arbitrary probability distributions, the Wasserstein distance between two Gaussian distributions admits a closed form expression~\cite{p:givens1984class}. In particular, given two multivariate Gaussian distributions $\mathcal{N}_1(\mu_1, \Sigma_1)$ and $\mathcal{N}_2(\mu_2, \Sigma_2)$, with $\mu_1$, $\mu_2\in\mathbb{R}^n$ and $\Sigma_1$, 
$\Sigma_2 \in \mathbb{S}^{++}_n$, the squared Wasserstein distance between them is given as follows:
\begin{align}\label{eq:wasserstein-definition}
& W^{2}\left(\mathcal{N}_1,\mathcal{N}_2\right)=\left\|\mu_{1}-\mu_{2}\right\|_{2}^{2} \nonumber \\
& \qquad \qquad+\operatorname{tr}\left(\Sigma_{1}+\Sigma_{2}-2\left(\Sigma_{2}^{1 / 2} \Sigma_{1} \Sigma_{2}^{1 / 2}\right)^{1 / 2}\right).
\end{align}
For more details, the reader may refer to \cite{p:halderacc16,p:balci2020}.

\subsection{Problem Setup}
We consider the following discrete-time stochastic linear system
\begin{subequations}
\begin{align}\label{eq:motion}
x(t+1) & = A(t) x(t) + B(t) u(t) + w(t), \\
x(0) & =x_0,~~ \quad ~~x_0 \sim \cN(\mu_0, \Sigma_0), \label{eq:motion2}
\end{align}
\end{subequations}
for $t \in [0,T-1]_{d}$, where $\mu_0\in
\mathbb{R}^n$ and $\Sigma_0 \in \mathbb{S}_n^{++}$ are given. Let 
 $\mathscr{X}_{0:t} := \{ x(\tau) \in \mathbb{R}^n: \tau \in [0,t]_{d} \}$, for $t\in [0,T]_d$, $\mathscr{U}_{0:t} := \{ u(\tau) \in
\mathbb{R}^m:~\tau\in [0, t]_{d}\}$, for $t\in [0,T-1]_d$, and $\mathscr{W}_{0:t} := \{ w(\tau) \in \mathbb{R}^n:~\tau\in
[0,t]_{d}\}$, for $t\in [0,T-1]_d$. We assume that the noise process $\mathscr{W}_{0:t}$ corresponds to a sequence of independent and identically distributed normal random variables
with
\begin{align}\label{eq:w1}
\mathbb{E} \left[ w(t) \right] & = \mathbf{0}, \qquad \mathbb{E} \left[ w(t)w(t)\t \right] = \delta(t,\tau) W,
\end{align}
for all $t,\tau\in [0, T-1]_{d}$, where $W\in\mathbb{S}_n^{+}$  and
$\delta(t,\tau):=1$, when $t=\tau$, and $\delta(t,\tau):=0$,
otherwise. By assuming that $W \in \mathbb{S}_n^{+}$, we cover the case in which $w(t) = D v(t)$ where $D \in \mathbb{R}^{n \times p}$ and $\mathrm{cov}_v(t, \tau) = \delta(t,\tau) V$ with $V \in\mathbb{S}_n^{++}$ because in the latter case, we would have $W = D V D\t \in \mathbb{S}_n^+$; thus there is no loss of generality in assuming that the noise vector has the same dimension with the state vector. Furthermore, $x_0$ is independent of $W_{0: T-1}$, that is,
\begin{align}\label{eq:x0w}
\mathbb{E}\left[ x_0  w(t)\t  \right] & =
\mathbf{0}, \qquad \mathbb{E}\left[ w(t) x_0\t  \right] = \mathbf{0},
\end{align}
for all $t \in [0, T-1]_d$. All random variables are defined on a (fixed) complete probability space $(\Omega, \mathfrak{F}, \mathbb{P})$. 

Equation~\eqref{eq:motion} can be written more compactly as follows:
\begin{align}\label{eq:bmx}
\bm{x} = \mathbf{G}_{\bm{u}} \bm{u} + \mathbf{G}_{\bm{w}} \bm{w} +
\mathbf{G}_{0} x_0,
\end{align}
where $\bm{x} := \mathrm{vertcat}(\mathscr{X}_{0:T})  \in \mathbb{R}^{(T+1)n}$, $\bm{u} := \mathrm{vertcat}(\mathscr{U}_{0:T-1}) \in\mathbb{R}^{Tm}$
and $\bm{w} := \mathrm{vertcat}(\mathscr{W}_{0:T-1}) \in\mathbb{R}^{Tn}$.
Furthermore,
\begin{align*}
& \mathbf{G}_{\bm{u}} :=  \begin{bmatrix} \mathbf{0} & \mathbf{0} & \dots & \mathbf{0} \\
B(0) &  \mathbf{0} & \dots & \mathbf{0} \\
\Phi(2,1)B(0) & B(1) & \dots & \mathbf{0} \\
\vdots & \vdots & \dots & \vdots \\
\Phi(T,1)B(0) & \Phi(T,2)B(1) &
\dots & B(T-1)\end{bmatrix} ,\\
&\mathbf{G}_{\bm{w}} := \begin{bmatrix} \mathbf{0} & \mathbf{0} & \dots & \mathbf{0} \\
I & \mathbf{0} & \dots & \mathbf{0} \\
\Phi(2,1) & I & \dots & \mathbf{0} \\
\vdots & \vdots & \dots & \vdots \\
\Phi(T,1) & \Phi(T,2) & \dots & I
\end{bmatrix},
~~\mathbf{G}_{0}:=
\begin{bmatrix} I \\ \Phi(1,0) \\ \vdots \\ \Phi(T,0)
\end{bmatrix},
\end{align*}
where $\Phi(t,\tau) := A(t-1)\dots A(\tau)$ and $\Phi(t,t) = I$,
for $t\in [1,T]_d$ and $\tau \in [0,t-1]_d$. Furthermore, we have
\begin{equation}\label{eq:bmw}
\mathbb{E}\left[ \bm{w} \right] = \bm{0},~~~\mathbb{E}\big[ \bm{w} \bm{w}\t \big] = \mathbf{W},
\end{equation}
where in light of \eqref{eq:w1}
\begin{align*}
 \mathbf{W} & := \mathbb{E}\big[
\mathrm{bdiag}(w(0)w(0)\t, \dots, w(T-1) w(T-1)\t \big] \\
& = \mathrm{bdiag}(W, \dots, W).
\end{align*}

\subsection{Affine Disturbance Feedback Controller Parametrization}
Under the assumption of perfect state information, one can recover at each stage the disturbance terms that have acted upon the system at all previous stages. Thus, one can use all these past distrurbances to compute the control input that will be applied to the system at each stage.
Next, we propose a modified version of the so-called affine disturbance feedback control policy parametrization, 
which we denote by ${\kappa}(t, W_{0:t-1}, x(0))$ 
which is defined as follows:
\begin{align}\label{eq:policy}
	\kappa(t)& = \begin{cases}
		\bar{u}(t) + L_{t} \left(x(0) - \mu_0\right) \\
		~ + \sum_{\tau = 0}^{t-1}K_{(t-1,\tau)}w(\tau)   & \mathrm{if}  ~ t \in [1, T-1]_{d}, \\
		\bar{u}(0) + L_{t} \left( x(0) - \mu_0 \right)  & \mathrm{if}  ~ t = 0,
	\end{cases}
\end{align}
where $L_{t}, K_{(t, \tau)} \in \mathbb{R}^{m \times n}$, $\forall t, \tau \in [0, T-1]_\mathrm{d}$. 

It is worth noting that the parametrization we propose is different than the standard affine disturbance feedback parametrization in that it has the extra term $L_t (x(0) - \mu_0)$ for every $t$. The reason why we added this term is to control the effect of the initial state covariance to the future state covariances. To better understand the necessity of this extra term, one can think of a disturbance-free system whose state covariance evolves as follows: $\mathrm{var}_{x(t+1)} = A(t) \mathrm{var}_{x(t)} A(t)\t$. In this case, the evolution of the state covariance  
cannot be affected by means of the standard affine disturbance feedback parametrization (there is no dependence on any control policy parameter). The situation is even worse for unstable systems in the presence of disturbances, where the disturbance feedback term $K_{(t, \tau)}$ would not be able to mitigate the growth of uncertainty.


If we set $u(t) = \kappa(t)$ for all $t \in [0, T-1]_d$ with $\kappa(t)$ as defined in \eqref{eq:policy}, then the resulting closed loop dynamics can be written as follows:
\begin{align}\label{eq:cloop}
x(t+1) & = A(t) x(t) + B(t) \bar{u}(t)  +  w(t)\nonumber \\
& B(t) \Big( L_{t}(x(0) - \mu_0) + \sum_{\tau =0}^{t-1} K_{(t-1,\tau)} w(\tau) \Big) 
\end{align}

\begin{remark}\label{remark:truncation}
To reduce the computational burden, one can truncate the disturbance history feedback policy \eqref{eq:policy} up to a desired number i.e. to use only a desired number of past disturbances that have acted upon the system to compute the control input. 
We denote by $\gamma \in [0, T]_\mathrm{d}$ a parameter that determines the length of the truncated history of disturbances such that the term $\sum_{\tau = 0}^{t-1} K_{(t-1, \tau)} w(\tau)$ that appears in \eqref{eq:policy} is replaced by the term $\sum_{\tau = t-(1+\gamma)}^{t-1} K_{(t-1, \tau)} w(\tau)$. It follows that only the last $\gamma + 1$ disturbance terms can be used in the modified control policy. 
To solve the SCCS and HCCS problems based on the truncation, the optimization (matrix) variable $\bm{\cK}$ which is defined in \eqref{eq:mathcalK} will have to be revised by setting the blocks $K_{(t-1, \tau)}$ equal to $\mathbf{0}$ for all $\tau < t-(1+\gamma)$.
\end{remark}

\subsection{Problem Formulation}
Next, we provide the precise formulations for the two variants of the covariance steering problems based on the control policy parametrization which is described in \eqref{eq:policy}.

\begin{problem}[Hard Constrained Covariance Steering]\label{problem1}
Let $\mu_\mathrm{d}, \mu_0 \in \mathbb{R}^n$, $\Sigma_\mathrm{d} \in \mathbb{S}^{++}_{n}$, $\rho \in \mathbb{R}^{+}$ be given. Consider the system described by~\eqref{eq:cloop}. Then, find the
collection of matrix gains $\mathscr{K} := \{ K_{(t,\tau)}, ~ L_{t} \in \mathbb{R}^{m \times n}: (t,\tau) \in
[0,T-2]\times [0,T-2],~t \geq \tau \}$
and the sequence of vectors $\overline{\mathscr{U}} := \{\bar{u}(0),
\dots, \bar{u}(T-1) \}$ that minimize the following performance
index:
\begin{align}\label{eq:perfindex}
    J_1(\mathscr{U}, \mathscr{K}) := \mathbb{E}\bigg[  \sum_{t=0}^{T-1}  u(t)\t u(t) \bigg]
\end{align}
subject to the boundary condition on the terminal mean and covariance:
\begin{equation}\label{eq:termmean1}
\mu_{x(T)} = \mu_\mathrm{d}. 
\end{equation}
\begin{equation}\label{eq:termcov1}
    \mathrm{var}_{x(T)} \preceq \Sigma_\mathrm{d}
\end{equation}
\end{problem}
\begin{remark}
The objective function defined in \eqref{eq:perfindex} represents the expected value of the total control input effort.
The constraint in \eqref{eq:termmean1} dictates that terminal state mean is equal to the desired mean and the constraint in \eqref{eq:termcov1} dictates that terminal state covariance is less then or equal to the desired one in the positive definite sense. If we want the terminal state covariance to precisely equal to the desired covariance, this formulation may not be adequate.\cite{p:BAKOLAS2018}
\end{remark}


\begin{problem}[Soft Constrained Covariance Steering]\label{problem2}
Let $\mu_\mathrm{f}, \mu_0 \in \mathbb{R}^n$, $\Sigma_\f, \Sigma_0 \in \mathbb{S}_{n}^{++}$ and $\rho \in \mathbb{R}^{+}$ be given. Consider the system described by~\eqref{eq:cloop}. Then, find the
collection of matrix gains $\mathscr{K} := \{ K_{(t,\tau)}, ~ L_{t} \in \mathbb{R}^{m \times n}: (t,\tau) \in
[0,T-2]\times [0,T-2],~t \geq \tau \}$
and the sequence of vectors $\overline{\mathscr{U}} := \{\bar{u}(0),
\dots, \bar{u}(T-1) \}$ that minimize the following performance
index:
\begin{equation}\label{eq:wassersteinobjective}
    J_{2}(\overline{\mathscr{U}}, \mathscr{K}) := W_{2}^{2}(\mathcal{N}_\mathrm{f}, \mathcal{N}_{\mathrm{d}}),
\end{equation}
subject to input constraint $C_{total}(\overline{\mathscr{U}}, \mathscr{K}) \leq 0$
, where
\begin{equation}\label{eq:totalinputeffort}
C(\overline{\mathscr{U}}, \mathscr{K}) := \mathbb{E}\Big[  \sum_{t=0}^{T-1}
u(t)\t u(t) \Big] - \rho^2,
\end{equation}
$\mathcal{N}_\mathrm{f} = \mathcal{N}(\mu_{x(T)}, \mathrm{var}_{x(T)})$ and $\mathcal{N}_\mathrm{d} = \mathcal{N}(\mu_{d}, \Sigma_{d})$ represent the Gaussian probability distribution of the terminal state at $t=T$ and the desired (goal) Gaussian probability distribution, respectively.
\end{problem}

\begin{remark}
Problem \ref{problem2} differs from the Problem \ref{problem1} in that the boundary conditions that are given in equations \eqref{eq:termmean1}, \eqref{eq:termcov1} were removed and replaced by a terminal cost term \eqref{eq:wassersteinobjective} and additional input constraints.
\end{remark}

\section{Reduction of the HCCS Problem into a Semidefinite Program}\label{s:reduction2convex}
To reduce Problems \ref{problem1} and \ref{problem2} into tractable optimization problems, first we need to express the control input vector $\bm{u}$ in terms of the decision variables defined in \eqref{eq:bmx}. In particular, the concatenated control input vector can be expressed as follows:
\begin{equation}\label{eq:bmu}
    \bm{u} = \bm{\Bar{u}} + \bm{\mathcal{L}} (x_0 - \mu_0) + \bm{\mathcal{K}} \bm{w}
\end{equation}
where
\begin{equation}
    \bm{\Bar{u}} := [\Bar{u}(0)\t, \Bar{u}(1)\t \dots, \Bar{u}(T-1)\t]\t,
\end{equation}
\begin{align}\label{eq:mathcalK}
	  \bm{K} := \begin{bmatrix} 
		K_{(0, 0)} & \mathbf{0}  &  \dots & \mathbf{0}  \\
		K_{(1, 0)} & K_{(1, 1)}  &  \dots & \mathbf{0}  \\
		\vdots & \vdots & \ddots & \vdots \\
		K_{(T-2, 0)} & K_{(T-2, 1)} & \dots & K_{(T-2, T-2)}
	\end{bmatrix},
\end{align}
\begin{equation}\label{eq:mathcalL}
    \bm{\mathcal{L}} := 
	\begin{bmatrix}
	L(0) \\
	L(1) \\
	\vdots \\
	L(T-1) 
	\end{bmatrix},
\end{equation} 
and $\bm{\mathcal{K}} = \left[ \begin{smallmatrix} \bm{0} & \bm{0} \\ \bm{K} & \bm{0} \end{smallmatrix} \right]$.

By plugging equation \eqref{eq:bmu} into equation \eqref{eq:bmx}, we get
\begin{align}\label{eq:compactbmx}
    &\bm{x} = \mathbf{G_u}\bm{\Bar{u}} + (\mathbf{G_w} + \mathbf{G_u} \bm{\mathcal{K}})\bm{w} + \mathbf{G_0}x_0 \nonumber \\ 
    & \qquad \qquad \qquad \qquad \qquad+ \mathbf{G_u} \bm{\mathcal{L}} (x(0) - \mu_0). 
\end{align}
Now, $x(t)$ can be expressed as $x(t) = \mathbf{P}_{t+1} \bm{x}$, where $\mathbf{P}_{t+1} := [\bm{0}, \dots , I_{n}, \dots, \bm{0}]$ is a block row vector whose $t^{th}$ block is equal to $I_n$ whereas all other blocks are equal to the zero matrix.

Next we provide analytical expressions for the mean and the variance of $\bm{x}$ and the state $x(t)$ for all $t \in [0, T-1]_d$.

\begin{proposition}\label{prop:meanvarx}
The mean and the variance of the random vector $\bm{x}$ which satisfies equation \eqref{eq:compactbmx} is given by:
\begin{equation}
    \mu_{\bm{x}} = \mathfrak{f}(\bm{\Bar{u}}), \qquad \mathrm{var}_{\bm{x}} = \mathfrak{h} (\bm{\mathcal{L}}, \bm{\mathcal{K}}),
\end{equation}
where 
\begin{subequations}
\begin{equation}\label{eq:fubar}
    \mathfrak{f}(\bm{\Bar{u}}) := \mathbf{G_u} \bm{\Bar{u}} + \mathbf{G_{0}} \mu_{0},
\end{equation}
\begin{align}
    &\mathfrak{h}(\bm{\mathcal{L}}, \bm{\mathcal{K}}) := (\mathbf{G_0} + \mathbf{G_u}\bm{\mathcal{L}}) \Sigma_0 (\mathbf{G_0} + \mathbf{G_u} \bm{\mathcal{L}})\t \nonumber \\
    & \qquad ~~~\qquad + (\mathbf{G_w} + \mathbf{G_u}\bm{\mathcal{K}}) \mathbf{W} (\mathbf{G_w} + \mathbf{G_u} \bm{\mathcal{K}})\t.
\end{align}
\end{subequations}
Furthermore, the mean and the variance of the state $x(t)$ are given by
\begin{subequations}
\begin{align}
     \mu_{x(t)} &= \mathbf{P}_{t+1} \mathfrak{f}(\bm{\Bar{u}}) \label{eq:termean} \\
     \mathrm{var}_{x(t)} &= \mathbf{P}_{t+1} \mathfrak{h}(\bm{\mathcal{L}}, \bm{\mathcal{K}}) \mathbf{P}_{t+1}\t. \label{eq:tervar}
\end{align}
\end{subequations}
\end{proposition}
The proof of this result and the proofs of other results in this paper are presented in the Appendix.

Next, we obtain an expression for the performance index of the minimum variance steering problem (Problem~\ref{problem1}) in terms of the decision variables $(\bar{\bm{u}}, \bm{\cL}, \bm{\cK})$.

\begin{proposition}\label{prop:J1}
The performance index $J_1(\overline{\mathscr{U}}, \mathscr{K})$ which is defined in \eqref{eq:perfindex} is equal to 
$\cJ_1(\bm{\Bar{u}}, \bm{\cL}, \bm{\cK})$, where
\begin{align}\label{eq:J1uK}
& \cJ_{1}(\bm{\Bar{u}}, \bm{\cL}, \bm{\cK}) := \bm{\Bar{u}}\t \bm{\Bar{u}} + \operatorname{tr}(\bm{\cK} \mathbf{W} \bm{\cK}\t) + \operatorname{tr}(\bm{\cL} \Sigma_0 \bm{\cL}\t)
\end{align}
provided that the pairs of decision variables $(\overline{\mathscr{U}},\mathscr{K})$ and $(\bar{\bm{u}},\bm{\cL}, \bm{\cK})$ are related by \eqref{eq:mathcalK}. Furthermore, $\cJ_1(\bm{\Bar{u}}, \bm{\cL}, \bm{\cK})$ is a convex function. 
\end{proposition}

The next proposition shows that terminal covariance constraint \eqref{eq:termcov1} can be written as a positive semidefinite constraint.
\begin{proposition}\label{prop:sdp}
The positive semi-definite constraint $\Sigma_{\mathrm{d}} \succeq \mathrm{var}_{x(T)}$ is satisfied iff $\bm{\mathcal{V}}(\bm{\cL}, \bm{\cK}) \in \mathbb{S}_{n}^{+}$ where 
\begin{equation}\label{eq:24}
    \bm{\mathcal{V}}(\bm{\cL}, \bm{\cK}) :=
    \begin{bmatrix}
    \Sigma_{\mathrm{d}} & \bm{\zeta}(\bm{\cL}, \bm{\cK}) \\ \bm{\zeta}(\bm{\cL}, \bm{\cK})\t & I_n
    \end{bmatrix}
\end{equation}
and $\bm{\zeta}(\bm{\cL}, \bm{\cK})$ is defined as in \eqref{eq:defzeta}.
\end{proposition}

To show that the terminal covariance constraint in \eqref{eq:termcov1} can be written as positive semidefinite constraint \eqref{eq:24}, we define:
\begin{equation}\label{eq:defzeta}
    \bm{\zeta}(\bm{\cL}, \bm{\cK}):= 
    \mathbf{P}_{T+1}
	\begin{bmatrix}
		\left( \mathbf{G}_{0} + \mathbf{G}_{\bm{u}} \bm{\cL} \right) & \left( \mathbf{G}_{\bm{w}} + \mathbf{G}_{\bm{u}} \bm{\cK} \right)
	\end{bmatrix}
	\mathbf{R}
\end{equation}
where $\mathbf{R} \mathbf{R}\t = \left[\begin{smallmatrix} \Sigma_0 & \mathbf{0} \\ \mathbf{0} & \mathbf{W} \end{smallmatrix}\right]$ and $\mathbf{P}_{T+1} \mathfrak{h}(\bm{\cL, \cK}) \mathbf{P}_{T+1}\t = \bm{\zeta}(\bm{\cL}, \bm{\cK})\bm{\zeta}(\bm{\cL}, \bm{\cK})\t$.

\begin{theorem}
The following semi-definite program is equivalent to Problem \ref{problem1}:
\begin{subequations}\label{eq:optprob1}
\begin{align}
    &\underset{\Bar{\bm{u}}, \bm{\cL}, \bm{\cK}}{\operatorname{min}}  && \cJ_{1}(\Bar{\bm{u}}, \bm{\cL}, \bm{\cK})  \qquad \\
    & \text{subject to}  && \mathbf{P}_{T+1} \mathfrak{f}(\bm{\Bar{u}}) = \mu_d \qquad  \\
    & && \bm{\mathcal{V}}(\bm{\cL}, \bm{\cK}) \in \mathbb{S}_{n}^{+}
\end{align}
\end{subequations}
\end{theorem}

\begin{remark}
Theorem \ref{theorem2} is a direct consequence of Propositions \ref{prop:meanvarx}, \ref{prop:J1} and \ref{prop:sdp}. In view of this theorem, Problem \ref{problem1} reduces into a semi-definite program.
\end{remark}


\section{Reduction of the SCCS Problem into a Difference of Convex Functions Program}\label{s:reduction2dc}
In this section, we associate the SCCS problem with a difference of convex functions program (DCP). In order to do that, we utilize the  control policy parametrization in \eqref{eq:policy} and use the results from Section \ref{s:reduction2convex}.

By setting $\mu_1 = \mu_{x(T)}$ and $\Sigma_1=\mathrm{var}_{x(T)}$, where $\mu_{x(T)}$ and $\Sigma_1=\mathrm{var}_{x(T)}$ are defined in \eqref{eq:termean} and \eqref{eq:tervar} for $t=T$, respectively, and also $\mu_2 = \mu_\f$ and $\Sigma_2 = \Sigma_{\mathrm{d}}$ into the expression of the squared Wasserstein distance given in \eqref{eq:wasserstein-definition}, we obtain the following expression of the objective function in terms of the new decision variables:
\begin{align}\label{eq:wasserstein-objective}
    & \cJ_3(\bm{\Bar{u},\cL, \cK}) := \lVert \mathbf{P}_{T+1} \mathfrak{f}(\bm{\Bar{u}}) - \mu_\mathrm{d} \rVert_2^2 \nonumber \\ 
    & \qquad ~~ + \operatorname{tr} \left( \mathbf{P}_{T+1}\mathfrak{h}(\bm{\cL}, \bm{\cK})\mathbf{P}_{T+1}\t + \Sigma_{\mathrm{d}} \right) \nonumber \\
    &  \qquad ~~ -2 \operatorname{tr}\left( (\sqrt{\Sigma_{\mathrm{d}}} \mathbf{P}_{T+1}\mathfrak{h}(\bm{\cL}, \bm{\cK})\mathbf{P}_{T+1}\t \sqrt{\Sigma_{\mathrm{d}}} )^{1/2} \right).
\end{align}

To show that the function defined in \eqref{eq:wasserstein-objective} is a difference of two convex functions we can define the objective function as a function of $\bm{\Bar{u}}$ and  $\bm{\zeta}$, where $\bm{\zeta}(\bm{\cL}, \bm{\cK})$ is an affine function which is defined in \eqref{eq:defzeta}.
In view of \eqref{eq:wasserstein-objective} and \eqref{eq:defzeta}, we define the new objective function as follows:
\begin{align}\label{eq:wasserstein-zeta}
    & \Tilde{\cJ}_3(\bm{\Bar{u}}, \bm{\zeta}(\bm{\cL}, \bm{\cK})) := \lVert \mathbf{P}_{T+1} \mathfrak{f}(\bm{\Bar{u}}) - \mu_\mathrm{d} \rVert_{2}^{2} \nonumber\\
    &\qquad  + \lVert \bm{\zeta}(\bm{\cL}, \bm{\cK}) \rVert_{F}^{2} + \operatorname{tr}(\Sigma_\mathrm{d})  -2 \lVert \sqrt{\Sigma_\mathrm{d}}\bm{\zeta}(\bm{\cL}, \bm{\cK}) \rVert_{*}.
\end{align}

\begin{proposition}\label{prop:wassersteinobjective}
The performance index $J_{2}(\overline{\mathscr{U}}, \mathscr{K})$ that is defined in equation \eqref{eq:wassersteinobjective} is equivalent to  $\Tilde{\cJ}_3(\bm{\Bar{u}}, \bm{\zeta}(\bm{\cL}, \bm{\cK}))$ which is defined in \eqref{eq:wasserstein-zeta}. Also, the function defined in the equation \eqref{eq:wasserstein-zeta} is the difference of two convex functions in variables 
$\bm{\Bar{u}}, \bm{\cL}, \bm{\cK}$.
\end{proposition}

Recall that, Problem \ref{problem2} has additional constraints compared with Problem \ref{problem1}. In the next proposition, we express these constraints and show that they correspond to convex constraints in terms of the decision variables ($\Bar{\bm{u}}, \bm{\cL}, \bm{\cK}$).
\begin{proposition}\label{prop:constraintCtotal}
The constraint function $C_{total}(\Bar{\mathscr{U}}, \mathscr{K})$ which is defined in equation \eqref{eq:totalinputeffort} can be expressed in terms of the decision variables $(\Bar{\bm{u}}, \bm{\cL}, \bm{\cK})$ as
\begin{align}
    & \cC(\bm{\Bar{u}}, \bm{\cL}, \bm{\cK}) := \bm{\Bar{u}}\t \bm{\Bar{u}} + \operatorname{tr}(\bm{\cK W \cK}\t)  \nonumber \\ 
    & \qquad \qquad \qquad \qquad + \operatorname{tr}(\bm{\cL} \Sigma_0 \bm{\cL}\t) - \rho^{2}.
\end{align}
and the set of all $(\bm{\Bar{u}, \cL, \cK})$ that satisfy the constraint $\cC(\bm{\Bar{u}, \cL, \cK}) \leq 0$ defines a convex set.
\end{proposition}



The next theorem which is a direct consequence of Propositions \ref{prop:wassersteinobjective} and \ref{prop:constraintCtotal} will allow us to reduce Problem \ref{problem2} into a difference of convex functions program:
\begin{theorem}\label{theorem2}
The following optimization problem is equivalent to Problem~\ref{problem2} :
\begin{subequations}\label{eq:optprob2}
\begin{align}
    && \underset{\bm{\Bar{u}}, \bm{\cL}, \bm{\cK}, \bm{\xi}}{\operatorname{min}} && \Tilde{\cJ}(\bm{\Bar{u}}, \bm{\zeta}) \qquad \qquad \quad ~~\\
    && \text{subject to} && \bm{\xi} = \bm{\zeta}(\bm{\cL}, \bm{\cK})  \qquad \qquad \qquad  \\
    && && \bm{\cC}(\bm{\Bar{u}}, \bm{\cL}, \bm{\cK}) \leq 0 \qquad \quad
\end{align}
where $\bm{\zeta}(\bm{\cL}, \bm{\cK})$ is defined in equation \eqref{eq:defzeta}.
\end{subequations}
\end{theorem}

One can exploit the structure of the performance index $\mathcal{J}_{3}$ that is defined in equation \eqref{eq:wasserstein-zeta} for efficient computation, by using the convex-concave procedure (CCP) which is a heuristic that is guaranteed to converge at stationary points of a difference of convex functions program \cite{p:yuille2003concave}. 
In the CCP, difference of convex objective and constraint functions are convexified by linearizing the difference function around the solution of the previous iteration. Then, the convexified problem is solved using convex optimization techniques. The procedure is terminated after the difference in the optimum objective values between iterations become sufficiently small.
To use the CCP as computational scheme for the DCP defined in Theorem \ref{theorem2}, the derivative of the term $-2 \lVert \sqrt{\Sigma_\mathrm{d}} \bm{\zeta}(\bm{\cL}, \bm{\cK}) \rVert_{*}$ in equation \eqref{eq:wasserstein-zeta} is required. 
Since the nuclear norm is a non-smooth function, this derivative may not exist in general but in this particular case, the derivative has a closed form expression which is given by the following proposition.

\begin{proposition}\label{prop:derivative}
If $\bm{\zeta} \bm{\zeta}\t \succ  \bm{0}, $ holds  $\forall \bm{\cL}, \bm{\cK}$ then  the gradient of $\lVert \sqrt{\Sigma_{\mathrm{d}}} \bm{\zeta}(\bm{\cL}, \bm{\cK}) \rVert_{*}$ is well-defined and given by:
\begin{align}
    &\nabla_{\zeta}\lVert \sqrt{\Sigma_{\mathrm{d}}} \bm{\zeta}(\bm{\cL}, \bm{\cK}) \rVert_{*} := \nonumber\\ 
    &\qquad \qquad\sqrt{\Sigma_{\mathrm{d}}} ( \sqrt{\Sigma_{\mathrm{d}}} \bm{\zeta}\bm{\zeta}\t \sqrt{\Sigma_{\mathrm{d}}})^{-1/2} \sqrt{\Sigma_{\mathrm{d}}} \bm{\zeta}.
\end{align}
\end{proposition}

\begin{remark}
Since the objective function of the optimization problem given in~\eqref{eq:optprob2} corresponds to a difference of convex function program and the constraints determine a convex set, the CCP heuristic is guaranteed to converge to a stationary point which satisfies the first order necessary optimality conditions. 
\end{remark}


\section{Numerical Expeeriments}\label{s:experiments}

In this section, we present results obtained by numerical experiments. We discuss various aspects of the proposed control policy parametrization including comparisons of its performance with the state history feedback parametrization used in our previous work \cite{p:BAKOLAS2018}. 
All computations were performed on a laptop with 2.8  GHz  Intel  Core  i7-7700HQ CPU  and  16  GB  RAM. We used CVXPY \cite{p:diamond2016cvxpy} as the modeling software for convex optimization problems and MOSEK \cite{mosek2010mosek} as the solver. To solve the SCCS problem, we used the convex-concave procedure by utilizing Proposition \ref{prop:derivative} to convexify the DCP objective function \eqref{eq:wassersteinobjective}. The termination criteria was chosen as $|f_k - f_{k-1}| \leq \epsilon$ where $\epsilon = 10^{-3}$ where $f_k$ denotes the result of the optimization problem at $k^{th}$ iteration for both examples.

In our numerical simulations, we consider two examples. One is based on a randomly generated linear dynamical system and the other corresponds to linearized model of the longitudinal dynamics of an aircraft. 
The randomly generated system was used for comparison purposes whereas the second system was chosen in order to show the applicability of our method in a more practical real-world problem.

\subsection{Random Linear System}
The parameters of the random linear system are taken as follows:
\begin{align*}
    A(t)= 
    \begin{bmatrix} 
    1.1 & -0.07 \\ 
    0.23 & -0.87 
    \end{bmatrix},
    \hfill
    B(t)= 
    \begin{bmatrix}
    0.0 \\
    0.1
    \end{bmatrix}, 
    \hfill
    W=
    \begin{bmatrix}
    0.1 & 0.0 \\
    0.0 & 0.3
    \end{bmatrix}
\end{align*}
Also, $\mu_0 = [1.0,0.0]\t$, $\Sigma_0 = I_2$, $\mu_\mathrm{d}=[10.0,0.0]\t$, $\Sigma_{\mathrm{d}}=\left[ \begin{smallmatrix} 4.0 & -1.5 \\ -1.5 & 4.0 \end{smallmatrix}\right]$ and $T=50$.
In our simulations, we have truncated the disturbance history in order to decrease the number of decision variables of the optimization problem as explained in Remark~\ref{remark:truncation}.

Table \ref{tab:totaltable} presents comparison results between the truncated version of the affine disturbance feedback control policy and the state history feedback control policy \cite{p:BAKOLAS2018}. The first column shows the truncation parameter $\gamma$ used in the controller parametrization of the former policy. 
The second and third columns show the objective value at the computed estimate of the minimizer and the computation time of the HCCS problem whereas the last two columns show the corresponding results obtained for the SCCS problem. 
The last row of the table shows the results obtained by using the state history feedback control policy. 

Based on these results, we can claim that the full disturbance history feedback policy achieves the same value for both the HCCS and the SCCS problems while reducing the computational cost. Also, as we increase the truncation parameter $\gamma$ (and thus increase the length of the truncated sequence of disturbances used in the control policy), we observe that the optimal value does not decrease below a certain value whereas the computation time increases.

The results based on the experiments with the random linear system suggest that our proposed policy parametrization when the full disturbance history is used may be equivalent to the state history feedback policy given that the two policies achieve the same objective value. However, more research is needed to establish rigorously the validity of the previous claim.

\begin{table}[ht]
    \centering
    \captionsetup{justification=centering}
    \caption{Comparison between affine disturbance feedback policy with different truncation lengths and state history feedback policy parametrization in terms of performance and computation time}
    \begin{tabular}{|c|c|c||c|c|}
      \multicolumn{1}{c}{ } & \multicolumn{2}{c}{HCCS}  &  \multicolumn{2}{c}{SCCS} \\ 
      \hline
      $\gamma$ & Value & Time (s) & Value & Time (s) \\
      \hline
      0 & 2839.93 & 1.32 & 2.00 & 2.19 \\
      \hline
      1 & 2535.38 & 0.89 & 1.19 & 3.63\\
      \hline
      3 & 2383.88 & 1.34 & 0.74 & 6.17\\
      \hline
      5 & 2333.49 & 1.69 & 0.59 & 9.78\\
      \hline
      10 & 2288.79 & 3.29 & 0.46 & 15.80\\
      \hline
      20 & 2271.51 & 4.84 & 0.40 & 29.77\\
      \hline
      30 & 2269.66 & 6.32 & 0.39 & 36.94\\
      \hline
      40 & 2269.46 & 9.78 & 0.38 & 38.87\\
      \hline
      50 & 2269.44 & 9.61 & 0.38 & 37.38\\
      \hline
      \hline
      \cite{p:bakcdc16} & 2269.44 & 16.32 & 0.38 & 112.39 \\
      \hline
    \end{tabular}
    \label{tab:totaltable}
\end{table}

\begin{figure*}[t!]
    \centering
    \begin{subfigure}{0.3\linewidth}
        \centering
        \includegraphics[width=\linewidth]{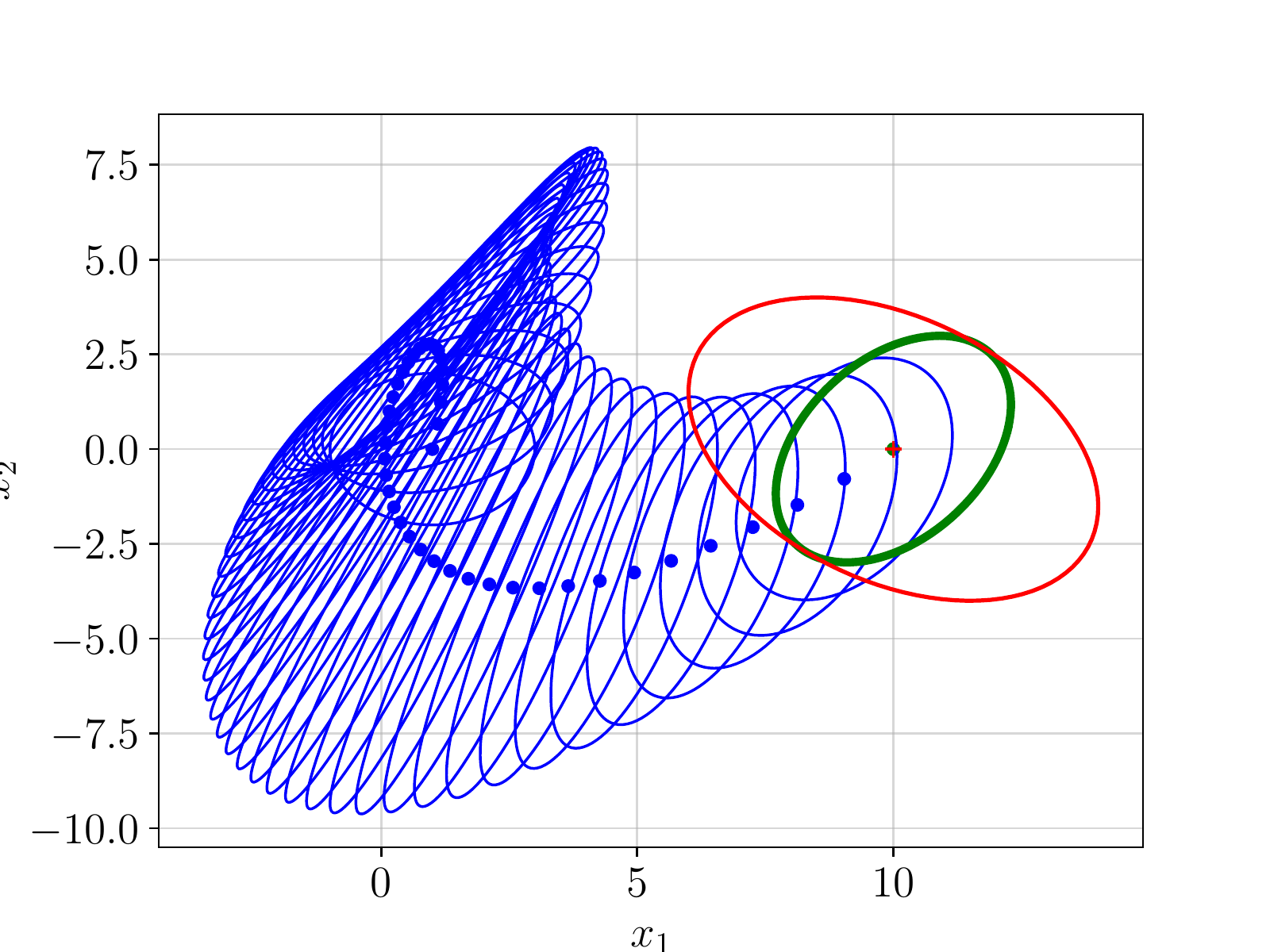}
        
    \end{subfigure}
    ~
    \begin{subfigure}{0.3\linewidth}
        \centering
        \includegraphics[width=\linewidth]{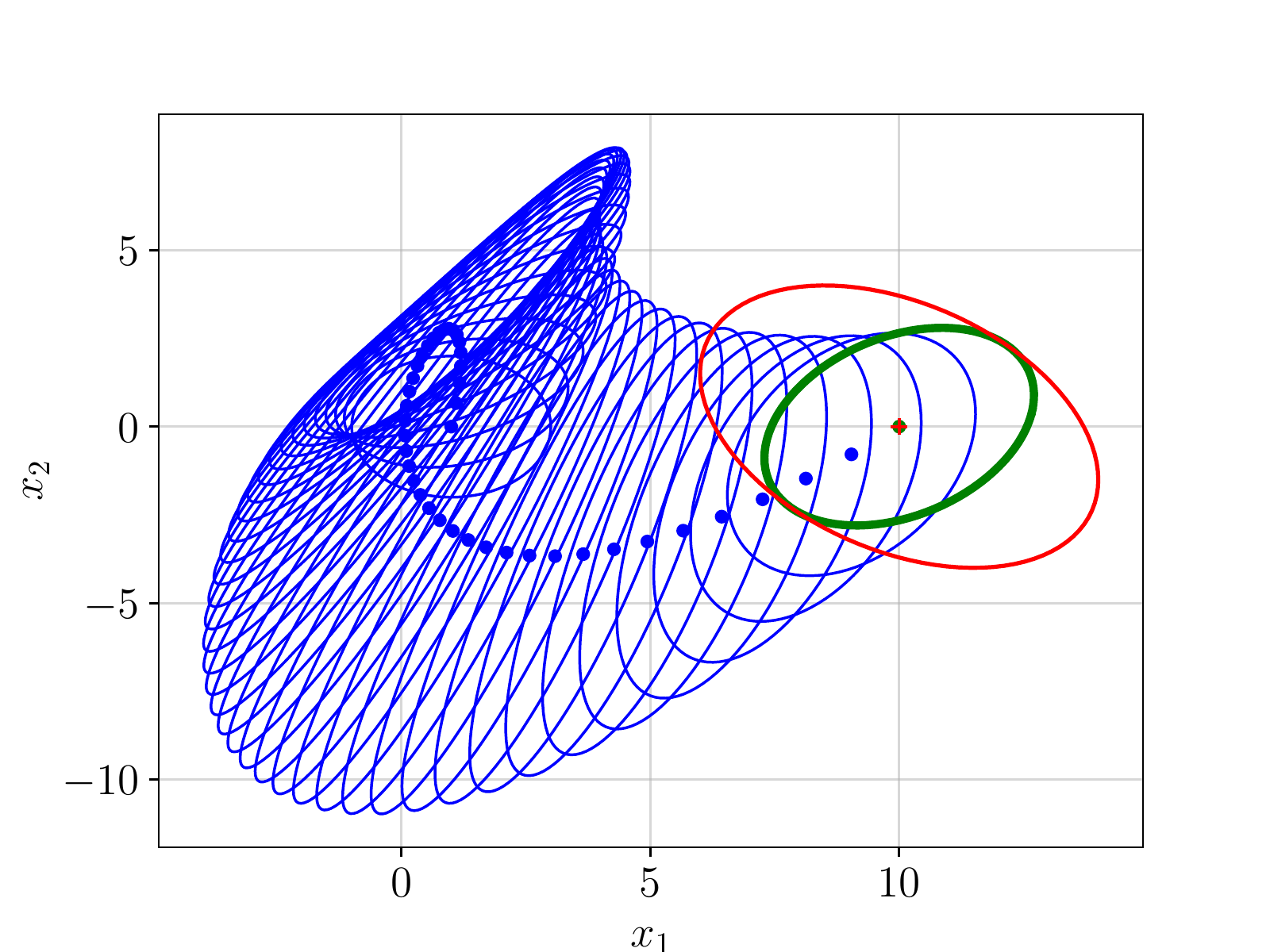}
        
    \end{subfigure}
    ~
    \begin{subfigure}{0.3\linewidth}
        \centering
        \includegraphics[width=\linewidth]{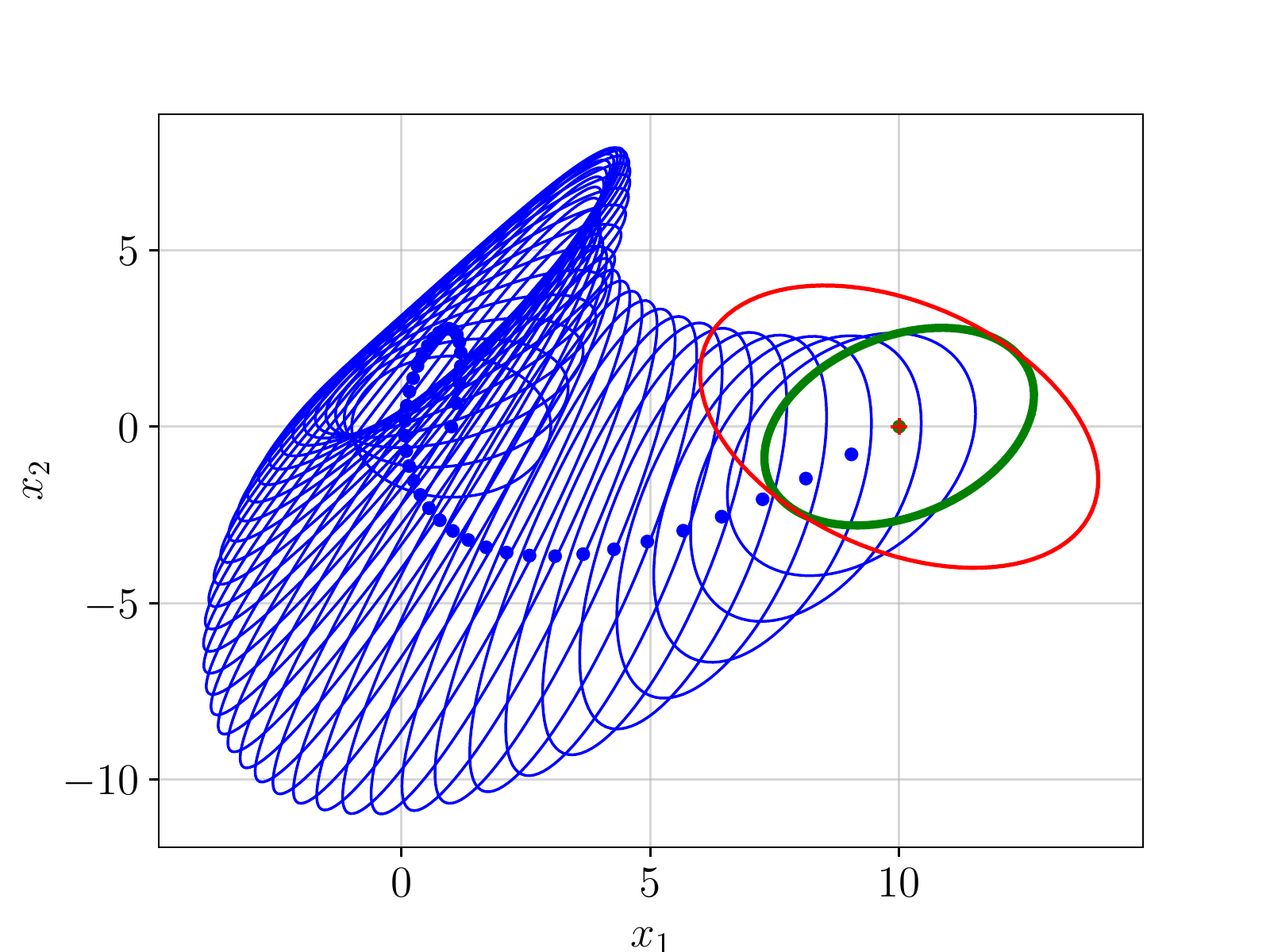}
        
    \end{subfigure}
    ~
    \begin{subfigure}{0.3\linewidth}
        \centering
        \includegraphics[width=\linewidth]{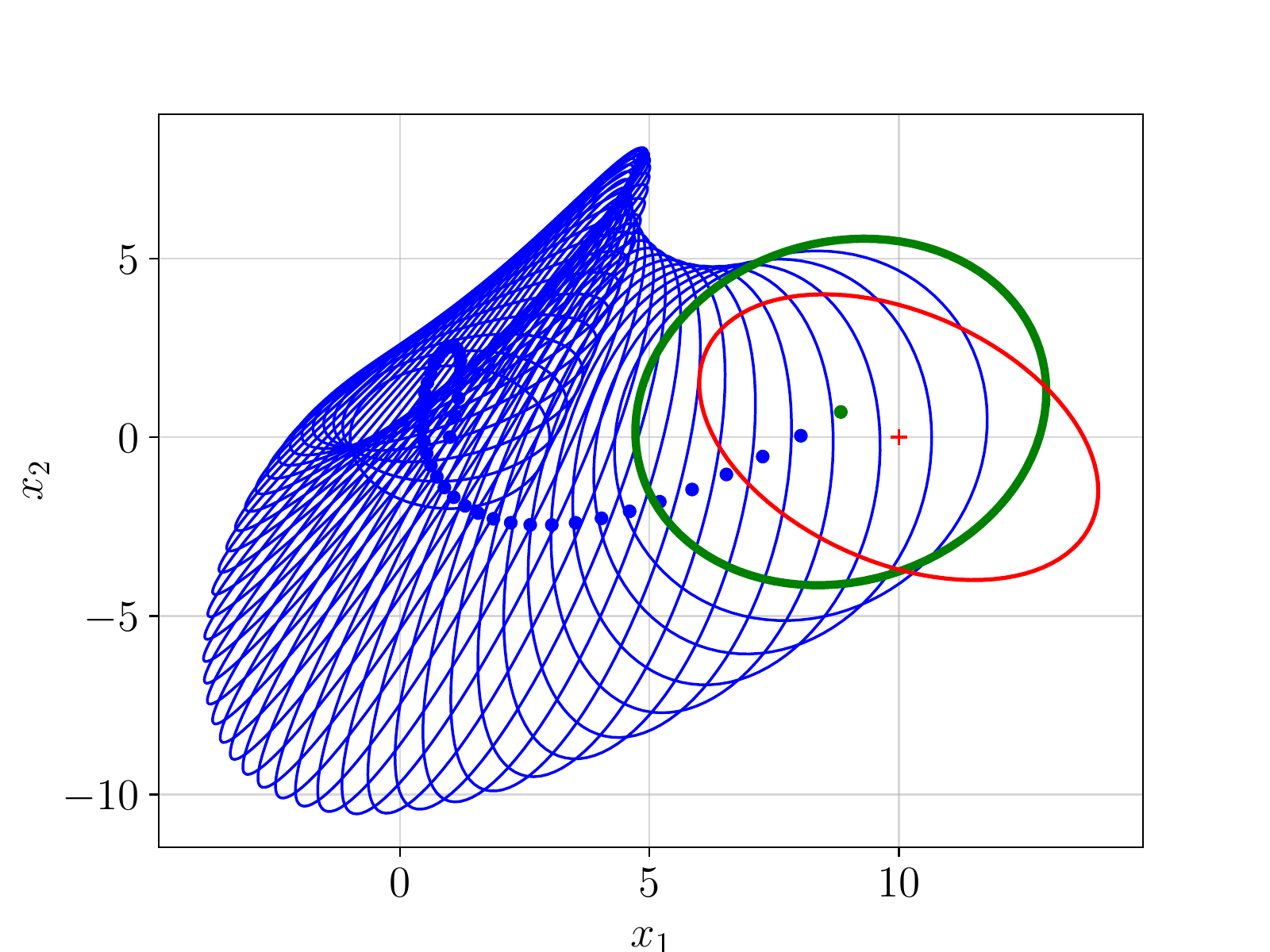}
        
    \end{subfigure}
    ~
    \begin{subfigure}{0.3\linewidth}
        \centering
        \includegraphics[width=\linewidth]{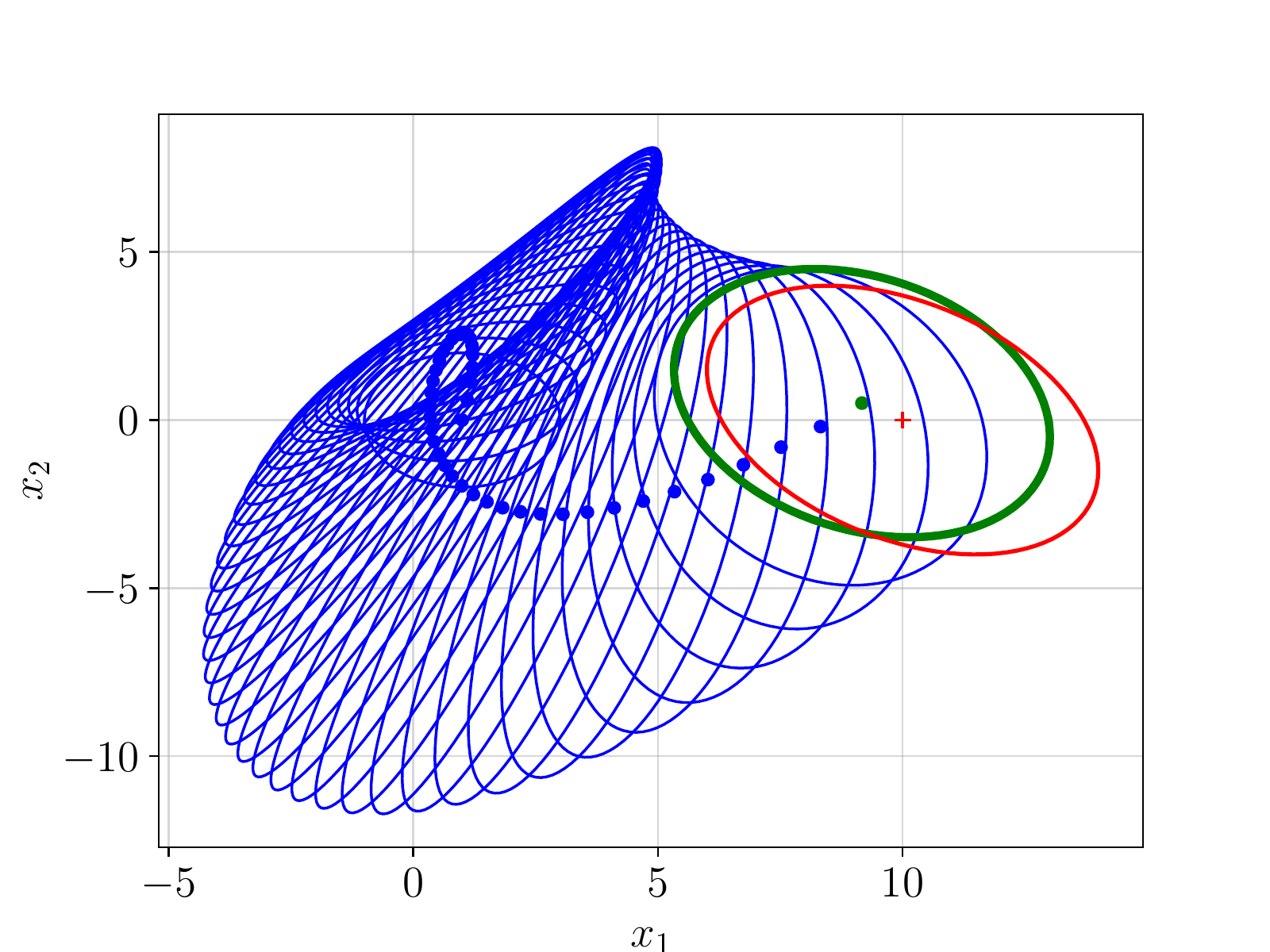}
        
    \end{subfigure}
    ~
    \begin{subfigure}{0.3\linewidth}
        \centering
        \includegraphics[width=\linewidth]{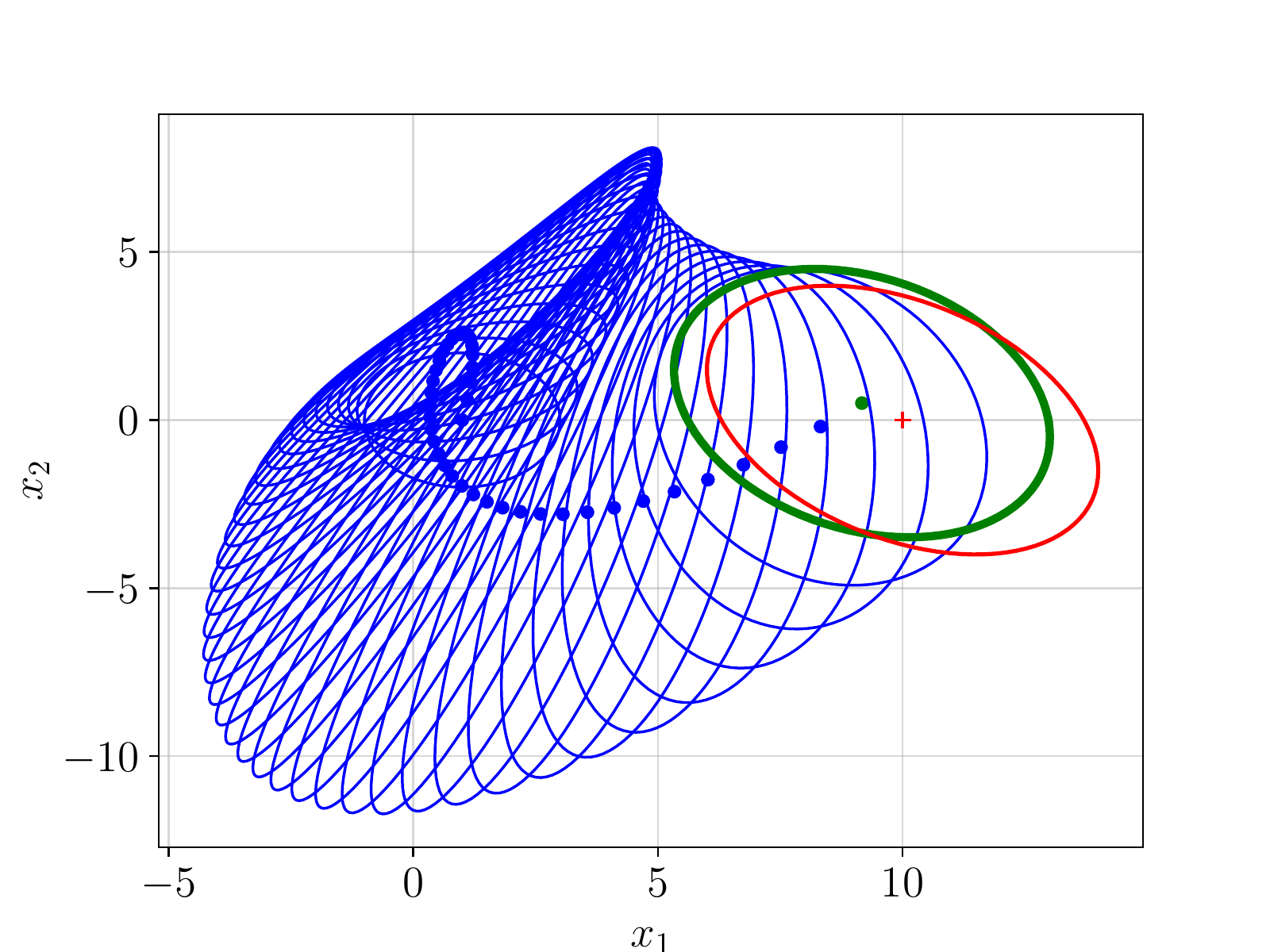}
        
    \end{subfigure}
    \caption{\small{Trajectories of the state mean and covariance of a random linear system: 
    The top plots illustrate the results for the HCCS problem and the bottom plots those for the SCCS problem. Only the last disturbance term is used in the control policy for the left subfigures ($\gamma=0$), the full disturbance history is used in the middle ones and the state history feedback policy is used on the right subfigures. Blue dots and ellipses show the evolution of the state mean and the 2-$\sigma$ confidence ellipse over time. The green dot and ellipse show the final distribution and the red cross and ellipse show the desired mean and the 2-$\sigma$ confidence ellipse of the desired terminal state distribution.} }
    \label{fig:traj-random}
\end{figure*}

\subsection{Linearized Longitudinal Aircraft Dynamics}
The discrete-time model of the linearized longitudinal dynamics of an aircraft is taken from \cite{p:kotsalis2020} and is obtained after the discretization of the continuous-time dynamics with a same sampling period ($\Delta T = 10$s). As initial and desired state mean we use $\mu_0 = [0.0, 0.0, 0.0, 0.0, 0.0]\t$ and $\mu_\mathrm{d} = [400.0, 0.0,0.0,0.0,0.0]\t$ whereas the initial and desired covariance matrices are chosen as $\Sigma_0 = \mathrm{bdiag}(100.0, 25.0, 25.0, 1.0, 1.0)$ and $\Sigma_\mathrm{d}= \mathrm{bdiag}(10^4, 100.0, 4.0, 1.0, 1.0)$. For the SCCS simulations, we take $\rho = 4$. The time horizon is chosen as $T = 40$.

\begin{figure}
    \centering
    \vspace{-0.5cm}
    \hspace{-1cm}
    \includegraphics[width=1.0\linewidth]{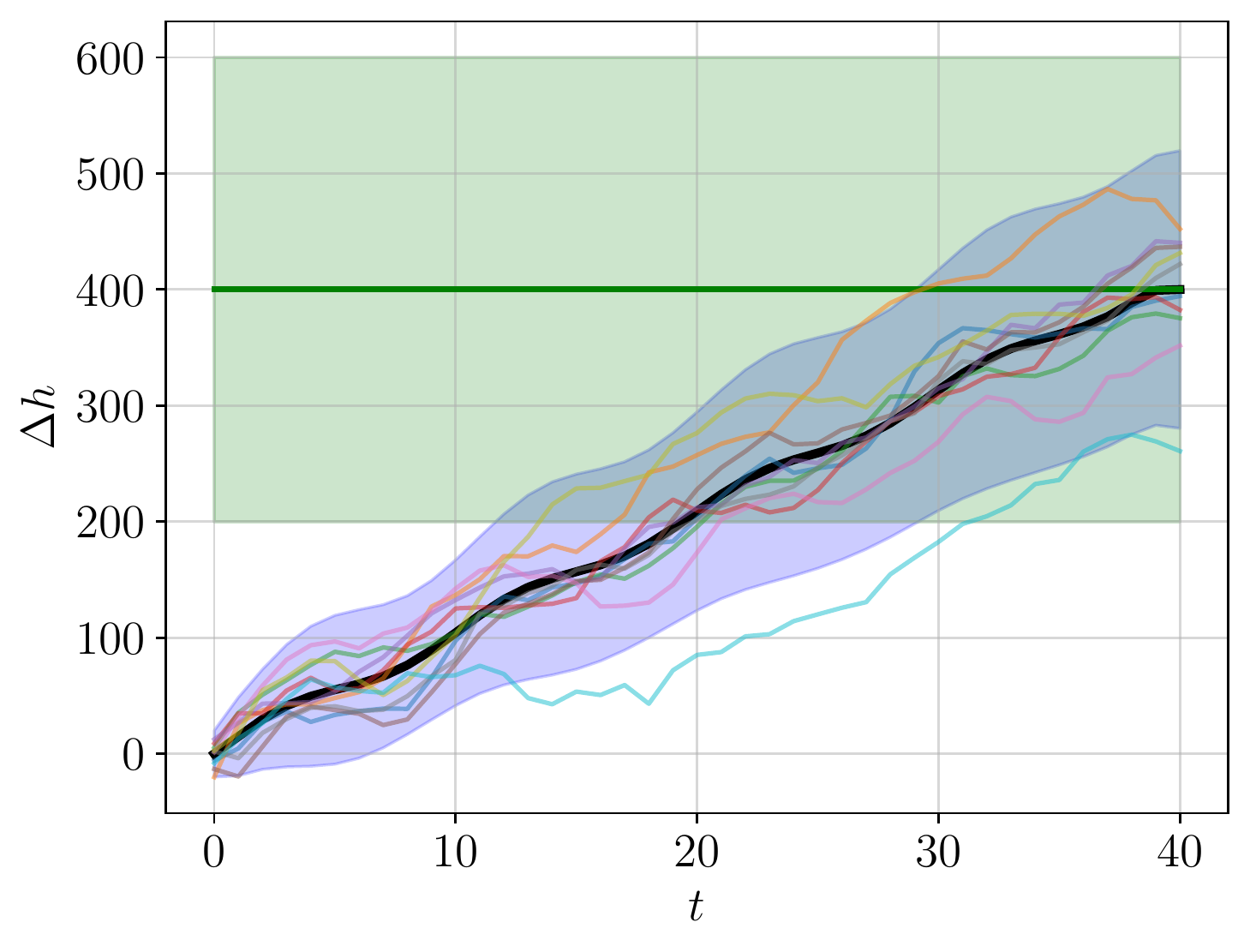}
    \caption{\small{$\Delta h$ vs $t$, (HCCS problem).}}
    \label{fig:aircraft-hard}
\end{figure}

\begin{figure}
    \centering
    \vspace{-0.5cm}
    \hspace{-1cm}
    \includegraphics[width=1.0\linewidth]{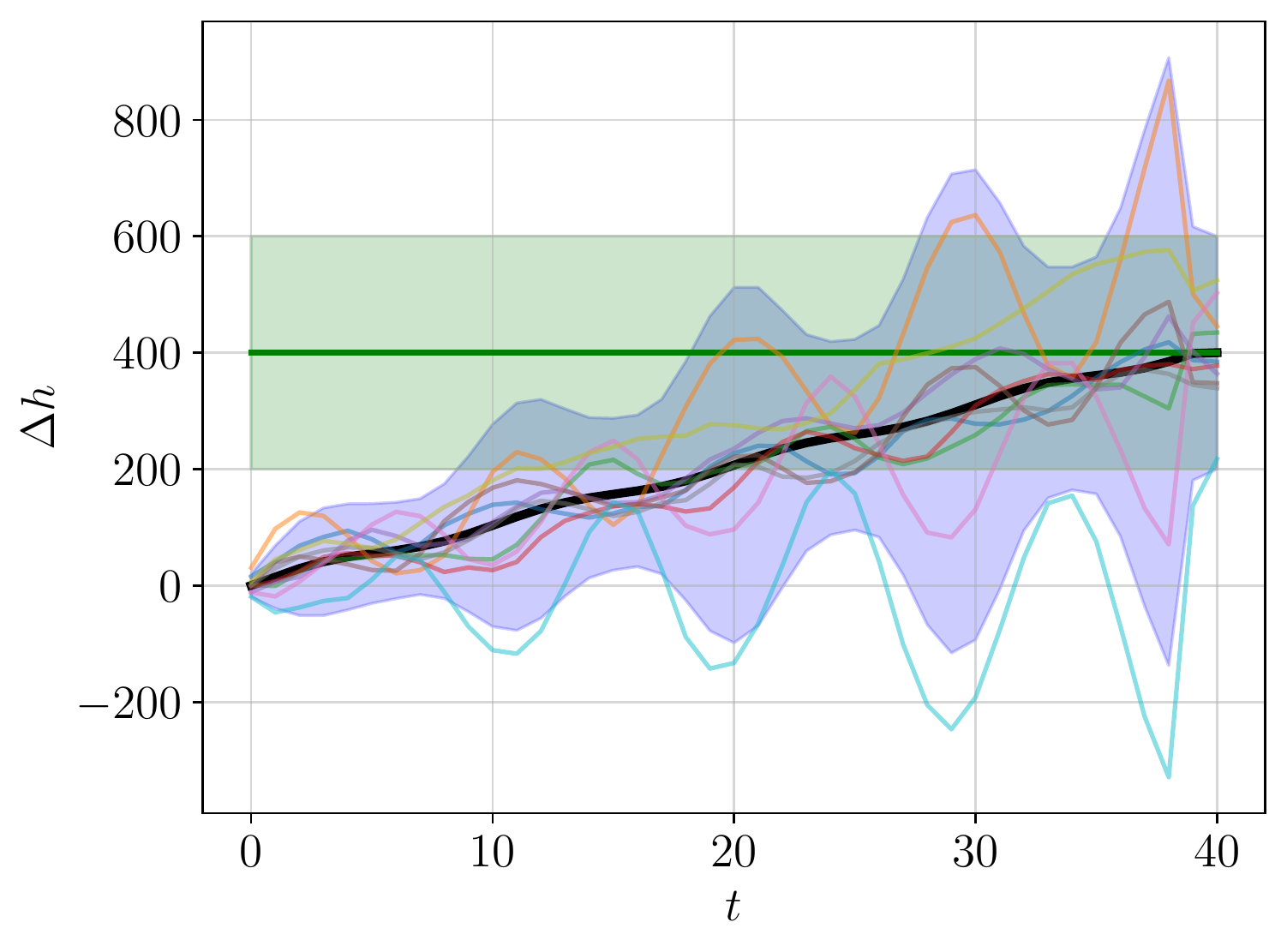}
    \caption{\small{$\Delta h$ vs $t$, (SCCS problem).}}
    \label{fig:aircraft-soft}
\end{figure}

In Figures \ref{fig:aircraft-hard} and \ref{fig:aircraft-soft}, the statistics of the first component of the state $x(t)$, which is the deviation from steady flight altitude denoted as $\Delta h$, is shown along with sample trajectories and the desired mean and covariance. 
Figure \ref{fig:aircraft-hard} corresponds to the HCCS problem whereas Figure~\ref{fig:aircraft-soft} to the SCCP problem. The black curve corresponds to the trajectory of the mean whereas the blue shaded area illustrates the 2-$\sigma$ confidence region, the green line and green shaded region shows the desired mean and the 2-$\sigma$ confidence region of the desired covariance respectively. The optimal policy is computed in 6.9s and 57.3s for the HCCS and the SCCS, respectively.
 

\section{Conclusion}\label{s:conclusion}
In this paper, we proposed a new covariance steering algorithm for discrete-time linear Gaussian stochastic systems. We have studied two variants of the covariance steering problem formulation to show the efficacy of the proposed parametrization. In both of these variants, the proposed parametrization performed better than the state history feedback policy parametrization in our numerical experiments. In our future work, we plan to expand this parametrization for distributed and robust covariance steering problems as well as covariance steering dynamic games.

\appendix

\subsection{Proof of Proposition \ref{prop:meanvarx}}
By applying the expectation operator on both sides of equation \eqref{eq:bmx}, we obtain:
\begin{align}\label{eq:Ebmx}
    \mu_{\bm{x}} &:= \mathbb{E}[\bm{x}] =\mathbf{G_u}\bm{\Bar{u}} + (\mathbf{G_w} + \mathbf{G_u} \bm{\mathcal{K}})\mathbb{E}[\bm{w}] + \mathbf{G_0}\mathbb{E}[x_0] \nonumber \\
    & \quad + \mathbf{G_u} \bm{\mathcal{L}} \mathbb{E}[(x(0) - \mu_0)] 
\end{align}
Since $\mathbb{E}[\bm{w}]=\bm{0}$, $\mathbb{E}[x_0] = \mu_0$, the equation \eqref{eq:fubar} follows immediately from equation \eqref{eq:Ebmx}. Furthermore, 
\begin{align}
    \bm{\Tilde{x}} &:= \bm{x} - \mu_{\bm{x}} \nonumber \\
    & = (\mathbf{G_w} + \mathbf{G_u}\bm{\cK})\mathbf{w} + (\mathbf{G_0} + \mathbf{G_u}\bm{\cL}) \Tilde{x}_0 \\
    \mathrm{var}_{\bm{x}} &= \mathbb{E}[\bm{\Tilde{x}}\bm{\Tilde{x}}\t] \nonumber \\
    & = \mathbb{E} \Big[ \left( \mathbf{G_{w}} + \mathbf{G_{u} \bm{\mathcal{K}}} \right) \bm{w} \bm{w}\t \left( \mathbf{G_{w}} + \mathbf{G_{u} \bm{\mathcal{K}}} \right)\t  \nonumber \\
	 & \quad\quad~ + \left( \mathbf{G_{w}} + \mathbf{G_{u} \bm{\mathcal{K}}} \right) \bm{w}  \Tilde{x}_0\t \left( \mathbf{G}_{0} + \mathbf{G_{u}} \bm{\mathcal{L}} \right)\t \nonumber \\
	 & \qquad~ + \left( \mathbf{G}_{0} + \mathbf{G_{u}} \bm{\mathcal{L}} \right) \Tilde{x}_0  \bm{w}\t \left( \mathbf{G_{w}} + \mathbf{G_{u} \bm{\mathcal{K}}} \right) \nonumber \\ 
	 & \qquad~ + \left( \mathbf{G}_{0} + \mathbf{G_{u}} \bm{\mathcal{L}} \right) \Tilde{x}_0 \Tilde{x}_0\t \left( \mathbf{G}_{0} + \mathbf{G_{u}} \bm{\mathcal{L}} \right)\t \Big].
\end{align}
Using the linearity of expectation functional and the identities $\mathbb{E}[\bm{w}] = 0$, $\mathbb{E}[\bm{w}\Tilde{x}_{0}\t] = 0$ and $\mathbb{E}[\Tilde{x}_0 \Tilde{x}_0\t] = \Sigma_0$, it follows:
\begin{align}
    \mathrm{var}_{\bm{x}} := & (\mathbf{G_0} + \mathbf{G_u} \bm{\cL}) \Sigma_0 (\mathbf{G_0} + \mathbf{G_u} \bm{\cL})\t \nonumber \\
    &~~ (\mathbf{G_w} + \mathbf{G_u} \bm{\cK}) \mathbf{W} (\mathbf{G_w} + \mathbf{G_u} \bm{\cK})\t.
\end{align}
Finally, using the expression $x(T) = \mathbf{P}_{T+1}\bm{x}$, we obtain:
\begin{align}
    \mathrm{var}_{x(T)} & = \mathbb{E}[(x(T)-\mu_{T})(x(T)-\mu_{T})\t] \nonumber \\
    & = \mathbf{P}_{T+1} \mathbb{E}[\bm{\Tilde{x}} \bm{\Tilde{x}}\t] \mathbf{P}_{T+1}\t \nonumber \\
    & = \mathbf{P}_{T+1} \mathfrak{h}(\bm{\cL}, \bm{\cK}) \mathbf{P}_{T+1}\t  \\
    \mu_{x(T)} & = \mathbf{P}_{T+1} \mathbb{E}[\bm{x}] = \mathbf{P}_{T+1} \mathfrak{f}(\bm{\Bar{u}}),
\end{align}
which proves the validity of Proposition \ref{prop:meanvarx}. $\hfill \blacksquare$

\subsection{Proof of Proposition \ref{prop:J1}}
We state that $\mathbb{E}[\sum_{0}^{T-1} u(t)\t u(t)] = \mathbb{E}[\bm{u}\t \bm{u}]$, since $\bm{u} = \mathrm{vertcat}(\mathscr{\overline{U}})$. Also, by using properties of the trace operator, it follows that:
\begin{align}
    \mathbb{E}[\bm{u}\t \bm{u}] &= \mathbb{E}[\operatorname{tr}(\bm{u}\t \bm{u})] = \operatorname{tr}(\mathbb{E} [\bm{u} \bm{u}\t]) \nonumber \\
    & = \operatorname{tr}(\mathbb{E}[ (\bar{\bm{u}} + \bm{\cK} \bm{w} + \bm{\cL} \tilde{x}_{0}) (\bar{\bm{u}} + \bm{\cK} \bm{w} + \bm{\cL} \tilde{x}_{0})\t ]). \nonumber
\end{align}
By using the identites $\mathbb{E}[\bm{w} \bm{w}\t] = \mathbf{W}$, $\mathbb{E}[\Tilde{x}_0 \bm{w}\t] = \bm{0}$ and $\mathbb{E}[\Tilde{x}_0 \tilde{x}_0 \t] = \Sigma_0$, we can further simplify the expression as:
\begin{align}
    \mathbb{E}[\sum_{t=0}^{T-1} u\t(t) u(t)] = \bm{\Bar{u}}\t\bm{\Bar{u}} + \operatorname{tr} (\bm{\cK} \mathbf{W} \bm{\cK}\t + \bm{\cL} \Sigma_0 \bm{\cL}\t)
\end{align}
which completes the proof. \hfill $\blacksquare$

\subsection{Proof of Proposition \ref{prop:sdp}}
In view of equations \eqref{eq:termean}, \eqref{eq:tervar}, the constraint $\mathrm{var}_x(T) \preceq \Sigma_\mathrm{d}$ can be written as the following positive semidefinite constraint $\cV(\bm{\cL}, \bm{\cK}) \in \mathbb{S}^{+}_n$ where,
\begin{align}\label{eq:covacons}
\cV(\bm{\cL}, \bm{\cK}) := \Sigma_\f - \mathbf{P}_{T+1} \mathfrak{h}(\bm{\cL}, \bm{\cK}) \mathbf{P}_{T+1}\t.
\end{align}
The expression of $\mathfrak{h}(\bm{\cL}, \bm{\cK})$ can alternatively be given by
\begin{align}\label{eq:huKalternate}
&\mathfrak{h}(\bm{\cL}, \bm{\cK}) = \nonumber\\
&~~ \begin{bmatrix}
	(\mathbf{G}_{0} + \mathbf{G}_{\bm{u}} \bm{\cL})\t \\
	(\mathbf{G}_{\bm{w}} + \mathbf{G}_{\bm{u}} \bm{\cK})\t
\end{bmatrix}\t
\begin{bmatrix}
	\Sigma_0 & \bm{0} \\
	\bm{0} & \mathbf{W}
\end{bmatrix}
\begin{bmatrix}
	(\mathbf{G}_{0} + \mathbf{G}_{\bm{u}} \bm{\cL})\t \\
	(\mathbf{G}_{\bm{w}} + \mathbf{G}_{\bm{u}} \bm{\cK})\t
\end{bmatrix}.
\end{align}
After plugging 
%
\eqref{eq:huKalternate} into \eqref{eq:covacons}, it follows that
\begin{equation}
 \cV(\bm{\cL}, \bm{\cK}) = \Sigma_\f - \bm{\zeta}(\bm{\cL}, \bm{\cK})\bm{\zeta}(\bm{\cL}, \bm{\cK})\t,
\end{equation}
where $\mathbf{R} \mathbf{R}\t := \left[\begin{smallmatrix}
    \Sigma_0 & \bm{0} \\
	\bm{0} & \mathbf{W}
\end{smallmatrix}\right] =: \mathbf{S}$ and 
\begin{subequations}
\begin{align}
	\bm{\zeta}(\bm{\cL}, \bm{\cK}) := 
	\mathbf{P}_{T+1}
	\begin{bmatrix}
		\left( \mathbf{G}_{0} + \mathbf{G}_{\bm{u}} \bm{\cL} \right) & \left( \mathbf{G}_{\bm{w}} + \mathbf{G}_{\bm{u}} \bm{\cK} \right)
	\end{bmatrix}
	\mathbf{R}. \nonumber
\end{align}
\end{subequations}                         
Thus, the constraint $\cV(\bm{\cL}, \bm{\cK}) \in \mathbb{S}^{+}_n$ can be written as follows:
\begin{equation*}
\mathbf{\cN}(\bm{\cL}, \bm{\cK})\in \mathbb{S}_{2n}^+,~~~\mathbf{\cN}(\bm{\cL}, \bm{\cK}) := \begin{bmatrix} \Sigma_\f & \bm{\zeta}(\bm{\cL}, \bm{\cK}) \\ \bm{\zeta}(\bm{\cL}, \bm{\cK})\t & I_n \end{bmatrix}.
\end{equation*}
The equivalence of the constraints $\cV(\bm{\cL}, \bm{\cK}) \in \mathbb{S}^{+}_n$ and $\mathbf{\cN}(\bm{\cL}, \bm{\cK}) \in \mathbb{S}_{2n}^+$ follows from the facts that $\Sigma_\mathrm{d} - \bm{\zeta (\cL, \cK)} I_n \bm{\zeta(\cL, \cK)} \t \in \mathbb{S}_{n}^{+}$, $I_n \in \mathbb{S}_{n}^{+}$ and the matrix $\cV(\bm{\cL, \cK})$ is the Schur complement of $I_n$ in $\mathcal{N}(\bm{\cL, \cK})$. Furthermore, given the function $\bm{\zeta}(\bm{\cL, \cK})$ is affine in ($\bm{\cL, \cK}$), thus the constraint $\mathcal{N}(\bm{\cL, \cK}) \in \mathbb{S}_{2n}^{+}$ is an LMI constraint. \hfill $\blacksquare$


\subsection{Proof of Proposition \ref{prop:wassersteinobjective}}
By using the equation \eqref{eq:termean} that is defined in Proposition \ref{prop:meanvarx}, we have a closed form expression for $\mu_{x(T)}$. Also, by using equations \eqref{eq:tervar} and \eqref{eq:defzeta} we have a closed form expression for $\mathrm{var}_{x(T)}$ as follows:
\begin{align}
    & \mu_{x(T)} = \mathbf{P}_{T+1} (\mathfrak{f}(\bm{\Bar{u}})),\qquad
    \mathrm{var}_{x(T)} =  \bm{\zeta}(\bm{\cL}, \bm{\cK}) \bm{\zeta} (\bm{\cL}, \bm{\cK}) \t. \nonumber
\end{align}
Plugging these equations into equation \eqref{eq:wasserstein-definition}, we have:
\begin{align}\label{eq:proofprop4}
    & \cJ_3(\bm{\Bar{u},\cL, \cK}) := \lVert \mathbf{P}_{T+1} \mathfrak{f}(\bm{\Bar{u}}) - \mu_\mathrm{d} \rVert_2^2 \nonumber \\ 
    & \qquad ~~ + \operatorname{tr} \left( \bm{\zeta}(\bm{\cL}, \bm{\cK}) \bm{\zeta}(\bm{\cL}, \bm{\cK})\t \right) + \operatorname{tr} \left(\Sigma_{\mathrm{d}} \right) \nonumber \\
    &  \qquad ~~ -2 \operatorname{tr}\Big( \left(\sqrt{\Sigma_{\mathrm{d}}} \bm{\zeta}(\bm{\cL}, \bm{\cK}) \bm{\zeta}(\bm{\cL}, \bm{\cK})\t \sqrt{\Sigma_{\mathrm{d}}} \right)^{1/2} \Big).
\end{align}
The second term in equation \eqref{eq:proofprop4} is equal to  $\lVert \bm{\zeta}(\bm{\cL}, \bm{\cK}) \rVert_{F}^{2}$ by the definition of the Frobenius norm. 
To show the equivalence of the forth term to the nuclear norm, let $\sqrt{\Sigma_{\mathrm{d}}} \bm{\zeta} = \mathbf{U S V}\t $ be the singular value decomposition of $\sqrt{\Sigma_{\mathrm{d}}} \bm{\zeta}$, thus
\begin{subequations}
\begin{align}
    & \sqrt{\Sigma_{\mathrm{d}}} \bm{\zeta} \bm{\zeta}\t  \sqrt{\Sigma_{\mathrm{d}}} = \mathbf{U} \mathbf{S S}\t \mathbf{U}\t \label{eq:proofprop4eq2}\\
    & \left( \sqrt{\Sigma_{\mathrm{d}}} \bm{\zeta} \bm{\zeta}\t  \sqrt{\Sigma_{\mathrm{d}}} \right)^{1/2} = \mathbf{U} \left( \mathbf{S S}\t \right)^{1/2} \mathbf{U}\t \label{eq:proofprop4eq3}\\
    & \operatorname{tr}\left( \left( \sqrt{\Sigma_{\mathrm{d}}} \bm{\zeta} \bm{\zeta}\t  \sqrt{\Sigma_{\mathrm{d}}} \right)^{1/2} \right) = \tr{\left(\mathbf{S S}\t \right)^{1/2}}.  \label{eq:proofprop4eq4}
\end{align}
\end{subequations}
The equation \eqref{eq:proofprop4eq2} uses the identity $\mathbf{V V}\t = I$ and the symmetry of $\sqrt{\Sigma_{\mathrm{d}}}$. The equality in \eqref{eq:proofprop4eq3} is a direct result from the definition of unique matrix square root of symmetric matrices. Finally, the equality in \eqref{eq:proofprop4eq4} is obtained by utilizing the identities $\operatorname{tr}(\mathbf{A B})=\operatorname{tr}(\mathbf{B A})$ and $\mathbf{U U}\t = I$. Since the matrix $\mathbf{S S}\t$ is  diagonal and its diagonal elements are the squared singular values of $\sqrt{\Sigma_{\mathrm{d}}} \bm{\zeta}$, the trace of $(\mathbf{S S}\t)^{1/2}$ is equal to the sum of singular values which is, by definition, equal to the nuclear norm. 

Lastly, both the frobenius norm and the nuclear norm are valid norms that satisfy the triangular equality thus, they are convex functions. Also, both $\mathfrak{f}(\bm{\Bar{u}})$ and $\bm{\zeta}(\bm{\cL}, \bm{\cK})$ are affine in their decision variables. Thus, the function defined in \eqref{eq:wasserstein-objective} is a difference of convex function in $(\bm{\Bar{u}}, \bm{\cL}, \bm{\cK})$. \hfill $\blacksquare$

\subsection{Proof of Proposition \ref{prop:constraintCtotal}}
It is trivial to establish the equivalence of $\bm{\cC}(\bm{\Bar{u}}, \bm{\cL}, \bm{\cK}) $ and $C_{total}(\overline{\mathscr{U}}, \mathscr{K})$ by using Proposition \ref{prop:J1} thus it is omitted. By using the properties of the Frobenius norm, we have  
\begin{align}\label{eq:proof5}
    &\bm{\cC}(\bm{\Bar{u}}, \bm{\cL}, \bm{\cK}) = \lVert \bm{\Bar{u}} \rVert_2^2 + \lVert \mathbf{R_w} \bm{\cK}\t \rVert_F^{2} \nonumber \\
    &\qquad\qquad\qquad\qquad\qquad 
    + \lVert \mathbf{R_0} \bm{\cL}\t \rVert_F^{2} - \rho^2
\end{align}
the first term on the right hand side of \eqref{eq:proof5} is a strictly convex quadratic function and it is thus convex. The second and third terms 
correspond to compositions of convex quadratic functions with affine functions and thus they are also convex. The last term is constant. We conclude that $\bm{\cC}(\bm{\Bar{u}}, \bm{\cL}, \bm{\cK})$ is a convex function. \hfill $\blacksquare$

\subsection{Proof of Proposition \ref{prop:derivative}}
To compute the derivative of $\lVert \sqrt{\Sigma_{\mathrm{d}}}\bm{\zeta} \rVert_{*}$, we use the matrix calculus tools given in \cite{b:magnus2019matrix}. For the sake of brevity, we denote $\sqrt{\Sigma_\mathrm{d}} \bm{\zeta \zeta}\t \sqrt{\Sigma_{\mathrm{d}}}$ by $X$ and $\lVert \sqrt{\Sigma_{\mathrm{d}}}\bm{\zeta} \rVert_{*}$ by $g(\bm{\zeta})$. We denote the differential operator by $\operatorname{d}$.
\begin{subequations}
\begin{align}
    g(\bm{\zeta}) = \tr{X^{1/2}}, \qquad  
    \d g(\bm{\zeta}) = \tr{\d X^{1/2}}. 
\end{align}
First take the differential of the identity $X^{1/2} X^{1/2} = X$:
\begin{align}
    \d X^{1/2} X^{1/2} + X^{1/2} \d X^{1/2}  = \d X.  \label{eq:propdereq1} 
\end{align}
By multiplying by $X^{-1/2}$ from the left both sides of the equation \eqref{eq:propdereq1}, it follows that
\begin{align}
    X^{-1/2} \d X^{1/2} X^{1/2} + \d X^{1/2} = X^{-1/2} \d X \\
    2 \tr{ \d X^{1/2}} = \tr{X^{-1/2} \d X}.
\end{align}
The expression for $X$ is quadratic in $\bm{\zeta}$, thus we take the differential of the expression of $\d X$ as follows:
\begin{align}
    \d X &= \d (\sqrt{\Sigma_{\mathrm{d}}} \bm{\zeta} \bm{\zeta}\t \sqrt{\Sigma_{\mathrm{d}}} ) \nonumber\\
    & = \sqrt{\Sigma_{\mathrm{d}}} \d \bm{\zeta} \bm{\zeta}\t \sqrt{\Sigma_\mathrm{d}} + \sqrt{\Sigma_{\mathrm{d}}} \bm{\zeta} \d \bm{\zeta}\t \sqrt{\Sigma_\mathrm{d}}. \label{eq:propdereq2}
\end{align}
Plugging the \eqref{eq:propdereq2} into \eqref{eq:propdereq1} we obtain the following:
\begin{align}
    \d g(\bm{\zeta}) &= \tr{\d X^{1/2}} \nonumber \\ 
    &= (1/2) 2 \tr{\bm{\zeta}\t \sqrt{\Sigma_\mathrm{d}} X^{-1/2} \sqrt{\Sigma_\mathrm{d}} \d \bm{\zeta}}
\end{align}
\end{subequations}
Finally, by using the Jacobian identification rule \cite{b:magnus2019matrix}, the derivative of $g(\bm{\zeta})$ is defined as:
\begin{equation}
    \nabla_{\bm{\zeta}} g(\bm{\zeta}) = \sqrt{\Sigma_\mathrm{d}} X^{-1/2} \sqrt{\Sigma_\mathrm{d}} \bm{\zeta}
\end{equation}
which completes the proof. \hfill $\blacksquare$

\bibliographystyle{ieeetr}
\bibliography{cdc2021}

\end{document}